\documentclass[12pt,oneside]{amsart}
\usepackage[utf8x]{inputenc}

\usepackage [a4paper, top=2.5cm, bottom=2.5cm, left=2cm, right=2cm, bindingoffset=8mm] {geometry}
\usepackage{amssymb,color,enumerate,amscd,anysize,graphicx}

\DeclareMathOperator{\R}{\mathbb{R}}
\DeclareMathOperator{\Z}{\mathbb{Z}}
\DeclareMathOperator{\cH}{\mathcal{H}}
\DeclareMathOperator{\Sp}{Sp}
\let\O\relax \DeclareMathOperator{\O}{O}
\DeclareMathOperator{\Sym}{Sym}
\DeclareMathOperator{\GL}{GL}
\DeclareMathOperator{\E}{\mathcal{E}}
\DeclareMathOperator{\tr}{tr}
\let\L\relax\DeclareMathOperator{\L}{\mathcal{L}}
\let\Im\relax\DeclareMathOperator{\Im}{Im}
\DeclareMathOperator{\SL}{SL}
\DeclareMathOperator{\e}{e}
\DeclareMathOperator{\SO}{SO}

\let\P\relax\DeclareMathOperator{\P}{\mathbb{P}}
\let\d\relax\DeclareMathOperator{\d}{d\!}

\newcommand{\lie}[1]{\mathfrak{#1}}
\let\t\relax \newcommand{\t}{^{t}\!}

\makeatletter
\renewcommand*\env@matrix[1][c]{\hskip -\arraycolsep
  \let\@ifnextchar\new@ifnextchar
  \array{*\c@MaxMatrixCols #1}}
\makeatother

\newtheorem{Thm}{Theorem}

\makeatletter                                                      
\def\tagform@#1{\maketag@@@{\ignorespaces#1\unskip\@@italiccorr}}
\makeatother

\title{Geometric classification of semidirect products \\ in the maximal parabolic subgroup of $\Sp(2,\R)$}
\author{F.~De~Mari}
\address{F. De Mari, DIMA\\
Via Dodecaneso, 35
\\16146 Genova, Italy}
\email{demari@dima.unige.it}

\author{E.~De Vito}
\address{E. De Vito, DIMA\\
Via Dodecaneso, 35
\\16146 Genova, Italy}
\email{devito@dima.unige.it}

\author{S.~Vigogna}
\address{S.~Vigogna, DIMA\\
Via Dodecaneso, 35
\\16146 Genova, Italy}
\email{vigogna@dima.unige.it.}
\date{July 18, 2013}

\input xy \xyoption{all}

\begin{document}

\begin{abstract}

 We classify up to conjugation by $\GL(2,\R)$ (more precisely, block diagonal symplectic matrices) all the semidirect products inside the maximal parabolic of $\Sp(2,\R)$ by means of an essentially  geometric argument. This classification has already been established in \cite{ABDD} without geometry, under a stricter notion of equivalence, namely conjugation by arbitrary symplectic matrices. The present approach might be useful in higher dimensions and provides some insight.

 \end{abstract}

 \maketitle
\pagestyle{myheadings}
\markboth{\sc Geometric classification}{\sc F. De Mari,  E. De Vito and S. Vigogna}

\section{Introduction}

 The main motivation for the classification problem we are after is the investigation of integral expansions of the type
 $$ \int_G \langle \varphi , \pi(g) \eta \rangle \pi(g) \eta  \d g,$$
 where $G$ is a locally compact group with  Haar measure $ \d g$ and where $\pi$ is a unitary representation of $G$
 on some Hilbert space $\cH$, typically $L^2(\R^d)$. The pair $(G,\pi)$ is said to be reproducing if the above
 superposition of one-dimensional projections is equal to $\varphi$ weakly. The primary examples of this are  wavelets and shearlets \cite{Dau92}, \cite{gukula06}.
Phase space analysis strongly indicates that when $\cH=L^2(\R^d)$, then, at least when geometric operations are to be encoded in the representation $\pi$, the correct setup is $G\subset\Sp(d,\R)$ and $\pi$ the restriction to $G$ of the metaplectic representation. We actually conjecture 
 that if such a pair is reproducing, then $G$ is conjugate (modulo compact factors) to a subgroup of 
 the (standard) maximal parabolic subgroup $Q$ of  the symplectic group, that is, the ``largest'' subgroup of block-lower triangular matrices. 
 This belief is supported by the many existing and very natural examples of this nature and, of course, by no known counterexamples. Particular attention deserves the so-called class $\E$, which consists of semidirect products $\Sigma\rtimes H$ inside $Q$, where $\Sigma$ is a vector part and $H$ is a linear subgroup of the general linear group. This class has been studied  by several authors in recent years 
 (see \cite{ABDD}, \cite{ADDM},  \cite{CT}, \cite{DD}),  and is a very rich source of examples. The classification, modulo conjugation by symplectic matrices, of the groups in $\E$ in the case $d=2$ has been carried out in the recent paper  \cite{ABDD} but the methods used there have no
 geometric flavour, as they are based mainly on Sylvester's law of inertia. Here we propose an intriguing  picture that might be of some help in the seemingly much harder problem of the analogous classification in $3D$. In order to keep the presentation clean and to convey the basic idea, the notion of equivalence adopted here is conjugation by block-diagonal matrices in $\Sp(2,\R)$, that is,  in the Langlands factor $MA$ of the standard maximal parabolic subgroup. This is the basic and most relevant conjugation, and was performed as a preliminary step in  \cite{ABDD}, where conjugation by any symplectic matrix is allowed, thereby giving rise to further identifications.

A general theory of  reproducing semidirect products of the same sort of those that we consider here has been developed by the first two listed authors in \cite{DD}. 

The main result here is Theorem~\ref{thm:classification}, whose proof is a consequence of what is established
in the subsequent sections.

\section{Preliminaries}
The symplectic group is 
\[
\Sp(d,\R)=\Bigl\{g\in\GL(2d,\R):\,^tgJg=J\Bigr\},
\]
where $J$ is the $2d\times2d$ matrix
\begin{equation} 
J = \begin{bmatrix}[r]
            0 & I \\
           -I & 0 
          \end{bmatrix}.
\label{J}
\end{equation}        
Its Lie  algebra  is
\[
\lie{sp}(d,\R)=\Bigl\{X\in\lie{gl}(2d,\R):\,^tXJ+JX=0\Bigr\}.
     \]
and its elements have the form
     \[
     \begin{bmatrix}
     A      &   \tau \\
     \sigma & - \t A
    \end{bmatrix},\;
     A \in \lie{gl}(d,\R), \quad \sigma,\tau \in \Sym(d,\R).
     \]
 The standard maximal parabolic subalgebra  $ \lie{q}$ of $\lie{sp}(d,\R)$ is  
\[
 \lie{q}=\Bigl\{
 \begin{bmatrix}
     A      &   0   \\
     x & - \t A
    \end{bmatrix}:
    A \in \lie{gl}(d,\R),\; x\in \Sym(d,\R)\Bigr\}
    \]
and is the Lie algebra of the maximal parabolic subgroup  $Q $ of $\Sp(d,\R)$ given by
\begin{equation}
Q=\Bigl\{
\begin{bmatrix}
     h        & 0        \\
     \sigma h & \t h^{-1}
    \end{bmatrix}:h \in \GL(d,\R),\; \sigma \in \Sym(d,\R)\Bigr\}.
\label{SH}
\end{equation}
Each $g\in Q$  can factorized (uniquely) as the product $g=\nu\mu\alpha$, namely
\begin{equation}
g= \begin{bmatrix}
     1      & 0 \\
     \sigma & 1
    \end{bmatrix} \begin{bmatrix}
                   m  & 0        \\
                   0  & \t m^{-1}
                  \end{bmatrix} \begin{bmatrix}
                                a & 0           \\
                                 0       &a^{-1}
                                \end{bmatrix}, \qquad 
                                  \det(m) = 1, \quad a\in \R,\,a\neq 0.
\label{langlands}\end{equation}
The groups to which these elements belong, namely
\begin{align*}
N&=\Bigl\{ \begin{bmatrix}
     1      & 0 \\
     \sigma & 1
    \end{bmatrix} : \sigma\in\Sym(d,\R)\Bigr\}, \\
M&=\Bigl\{ \begin{bmatrix}
                   m  & 0        \\
                   0  & \t m^{-1}
                  \end{bmatrix}: m\in\SL(d,\R)\Bigr\},    \\
A&=\Bigl\{  \begin{bmatrix}
                                a & 0           \\
                                 0       &a^{-1}
                                \end{bmatrix}: a\in\R\Bigr\}     \\
\end{align*}
are the factors of the  Langlands decomposition of  $Q$, which can be written either as  $Q = NMA$, like in \eqref{langlands}, or as $Q=MAN$. This means that each element factors as  $g=\nu\mu\alpha$ for unique elements $\nu\in N$, $\mu\in M$ and $\alpha \in A$, as in~\eqref{langlands},  or as $g=\mu'\alpha'\nu'$, with different elements but again unique with respect to this  order.
Notice that $ N \simeq \Sym(d,\R) $ and $ MA \simeq \GL(d,\R)\subset\Sp(d,\R)$, where the latter inclusion is given by the natural embedding $h\mapsto{\rm diag}(h,\t h^{-1})$. Parametrizing the elements of $Q$ as in \eqref{SH}, we  get
\[
(\sigma,h)(\sigma^\prime,h^\prime)=(\sigma+h^{\dag}(\sigma^\prime), hh^\prime)
\]
where 
 \begin{equation} \label{eq:dag_action}
h^{\dag}(\sigma) :=\ \t h^{-1}\sigma h^{-1}.
 \end{equation}
Thus, $Q$ is the semidirect product of the symmetric and the invertible matrices:
\[
 Q = \Sym(d,\R) \rtimes \GL(d,\R). 
 \]
The semidirect structure is generally not enjoyed by the subgroups of $Q$ (see \cite{ABDD} for an explicit example). As in \cite{ABDD}, \cite{ADDM},  \cite{CT},  we  define the class $\E_d$ as the family of all connected subgroups of $Q$ which are semidirect products, namely of the form
 $$ 
 \Sigma \rtimes H, \qquad 0 \neq \Sigma\subset\Sym(d,\R), \quad 1 \neq H \subset \GL(d,\R),
  $$
  where $\Sigma$ is a vector space and $H$ is a connected subgroup.

 While $N$ does not preserve $\E_d$, $MA$-conjugations, that is conjugations with matrices in the general linear group, do send groups in  
 $\E_d$ into groups in $\E_d$. More precisely, identifying again  $g$ with ${\rm diag}(g,\t g^{-1})$, we have
 \begin{equation} \label{eq:conjugation_split}
  g (\Sigma \rtimes H) g^{-1} = (\t g^{-1} \Sigma g^{-1}) \rtimes (g H g^{-1}).
 \end{equation}

The purpose of this article is to achieve a complete geometric
 classification of all the groups in $\E_2$ up to conjugation by
 elements in  $\GL(2,\R)\subset\Sp(2,\R)$.  The full list is provided
 by the following theorem, whose proof is the content of the remaining sections.

\begin{Thm}\label{thm:classification}
Up to conjugation by   $\GL(2,\R)$, any group $G$ in $\E_2$ is one of the  following:
 \begin{enumerate}[a)]
 \item  $5$-d groups
    \begin{align*}
    &  \Bigl\{ \begin{bmatrix}
        0 & u \\
        u & v
      \end{bmatrix} : u,v \in \R \Bigr\} \rtimes
      \Bigl\{ \begin{bmatrix}
        \e^t &  r \\
        0 & \e^s
      \end{bmatrix} : t,s,r \in \R \Bigr\}
    \end{align*}
 \item  $4$-d groups
   \begin{align*}
  &   \Bigl\{ \begin{bmatrix}
       u & v \\
       v & -u
     \end{bmatrix} : u,v \in \R \Bigr\} \rtimes \Bigl\{
     \e^t \begin{bmatrix}
       \cos\theta & -\sin\theta \\
       \sin\theta & \cos\theta
     \end{bmatrix} : t \in \R , \theta \in \mathbb{T} \Bigr\}
\\
& \Bigl\{ \begin{bmatrix}
       u & 0 \\
       0 & 0
     \end{bmatrix} : u \in \R \Bigr\} \rtimes \Bigl\{ \begin{bmatrix}
       \e^t &  0 \\
       r & \e^s
     \end{bmatrix} : t,s,r \in \R \Bigr\} 
\\
&     \Bigl\{ \begin{bmatrix}
       0 & u \\
       u & v
     \end{bmatrix} : u,v \in \R \Bigr\} \rtimes
     \Bigl\{ \begin{bmatrix}
       \e^{t\lambda} &  s \\
       0 & \e^{t(\lambda +1)}
     \end{bmatrix} :t,s \in \R\Bigr\} \quad \lambda \in \R 
\\
&     \Bigl\{ \begin{bmatrix}
       0 & u \\
       u & v
     \end{bmatrix} : u,v \in \R \Bigr\} \rtimes
     \Bigl\{ \begin{bmatrix}
       \e^t &  0  \\
       0 & \e^s
     \end{bmatrix}:t,s \in \R\Bigr\} 
\\ 
 &    \Bigl\{ \begin{bmatrix}
       0 & u \\
       u & v
     \end{bmatrix} : u,v \in \R \Bigr\} \rtimes
     \Bigl\{ \begin{bmatrix}
       \e^t &  s  \\
       0 & \e^t
     \end{bmatrix} :t,s \in \R\Bigr\} 
\\
 &  \Bigl\{ \begin{bmatrix}
       u & 0 \\
       0 & v
     \end{bmatrix} : u,v \in \R \Bigr\}  \rtimes
     \Bigl\{ \begin{bmatrix}
       \e^t &  0 \\
       0 & \e^s
     \end{bmatrix} t,s \in \R \Bigr\} 
   \end{align*}
\item $3$-d groups
\begin{align*}
& \Bigl\{ \begin{bmatrix}
      u & 0 \\
      0 & u
    \end{bmatrix} : u \in \R \Bigr\} \rtimes \Bigl\{
    \e^t \begin{bmatrix}
      \cos\theta & -\sin\theta \\
      \sin\theta & \cos\theta
    \end{bmatrix} : t \in \R , \theta \in \mathbb{T} \} 
\\
& \Bigl\{ \begin{bmatrix}
      u & v \\
      v & -u
    \end{bmatrix} : u,v \in \R \Bigr\} \rtimes \Bigl\{\e^{t
    } \begin{bmatrix}
      \cos t\alpha & -\sin t\alpha \\
      \sin t\alpha & \cos t\alpha
    \end{bmatrix} : t \in \R \Bigr\} \quad \alpha\geq0
\\
&    \Bigl\{ \begin{bmatrix}
      u & v \\
      v & -u
    \end{bmatrix} : u,v \in \R \Bigr\} \rtimes \Bigl\{ \begin{bmatrix}
      \cos\theta & -\sin\theta \\
      \sin\theta & \cos\theta
    \end{bmatrix} : \theta \in \mathbb{T} \Bigr\} 
\\
& \Bigl\{ \begin{bmatrix}
      u & 0 \\
      0 & 0
    \end{bmatrix} : u \in \R \Bigr\} \rtimes
    \Bigl\{ \begin{bmatrix}[l]
      \e^{t\lambda} &  0 \\
      s & \e^{t(\lambda +1)}
    \end{bmatrix} : t,s \in \R \Bigr\} \quad \lambda \in \R
\\
&    \Bigl\{ \begin{bmatrix}
      u & 0 \\
      0 & 0
    \end{bmatrix} : u \in \R \Bigr\} \rtimes \Bigl\{ \begin{bmatrix}
      \e^t &  0 \\
      0 & \e^s
    \end{bmatrix}:t,s \in \R\Bigr\} 
\\ 
&    \Bigl\{ \begin{bmatrix}
      u & 0 \\
      0 & 0
    \end{bmatrix} : u \in \R \Bigr\} \rtimes \Bigl\{ \begin{bmatrix}
      \e^t &  0 \\
      s & \e^t
    \end{bmatrix} :t,s \in \R\Bigr\} 
\\
&    \Bigl\{ \begin{bmatrix}
      0 & u \\
      u & v
    \end{bmatrix} : u,v \in \R \Bigr\} \rtimes \Bigl\{
    \e^{t\lambda} \begin{bmatrix}
      1 &  0  \\
      0 & \e^t
    \end{bmatrix}:t \in \R\Bigr\} \quad \lambda \in \R
\\
&    \Bigl\{ \begin{bmatrix}
      0 & u \\
      u & v
    \end{bmatrix} : u,v \in \R \Bigr\} \rtimes \Bigl\{ \begin{bmatrix}
      1 & t \\
      0 & 1
    \end{bmatrix}:t \in \R\Bigr\} 
\\
& \Bigl\{ \begin{bmatrix}
      0 & u \\
      u & v
    \end{bmatrix} : u,v \in \R \Bigr\} \rtimes \Bigl\{
    \e^t \begin{bmatrix}
      1 & t \\
      0 & 1
    \end{bmatrix}:t \in \R\Bigr\} 
\\
&    \Bigl\{ \begin{bmatrix}
      0 & u \\
      u & v
    \end{bmatrix} : u,v \in \R \Bigr\} \rtimes \Bigl\{ \begin{bmatrix}
      \e^t &  0 \\
      0 & \e^t
    \end{bmatrix} : t \in \R\Bigr\} 
\\
&    \Bigl\{ \begin{bmatrix}
      0 & u \\
      u & 0
    \end{bmatrix} : u \in \R \Bigr\} \rtimes \Bigl\{ \begin{bmatrix}
      \e^t &  0 \\
      0 & \e^s
    \end{bmatrix} t,s \in \R \Bigr\} 
\\
&    \Bigl\{ \begin{bmatrix}
      u & 0 \\
      0 & v
    \end{bmatrix} : u,v \in \R \Bigr\} \rtimes \Bigl\{ \begin{bmatrix}
      \e^t &  0     \\
      0 & \e^{t\beta}
    \end{bmatrix}:t \in \R\Bigr\}\quad \beta \in [-1,1]
  \end{align*}
\item $2$-d groups
  \begin{align*}
 &   \Bigl\{ \begin{bmatrix}
      u & 0 \\
      0 & u
    \end{bmatrix} : u \in \R \Bigr\} \rtimes \Bigl\{\e^{t
    } \begin{bmatrix}
      \cos t\alpha & -\sin t\alpha \\
      \sin t\alpha & \cos t\alpha
    \end{bmatrix} : t \in \R \Bigr\} \quad\alpha\geq0
\\
&    \Bigl\{ \begin{bmatrix}
      u & 0 \\
      0 & u
    \end{bmatrix} : u \in \R \Bigr\} \rtimes \Bigl\{ \begin{bmatrix}
      \cos\theta & -\sin\theta \\
      \sin\theta & \cos\theta
    \end{bmatrix} : t \in \R , \theta \in \mathbb{T} \Bigr\} 
\\
& \Bigl\{ \begin{bmatrix}
      u & 0 \\
      0 & 0
    \end{bmatrix} : u \in \R \Bigr\} \rtimes \Bigl\{
    \e^{t\lambda} \begin{bmatrix}
      1 &  0  \\
      0 & \e^t
    \end{bmatrix}:t \in \R\Bigr\} \quad \lambda \in \R 
\\
&\Bigl\{ \begin{bmatrix}
      u & 0 \\
      0 & 0
    \end{bmatrix} : u \in \R \Bigr\} \rtimes \Bigl\{ \begin{bmatrix}
      1 & 0 \\
      t & 1
    \end{bmatrix}:t \in \R\Bigr\} 
\\
&  \Bigl\{ \begin{bmatrix}
      u & 0 \\
      0 & 0
    \end{bmatrix} : u \in \R \Bigr\} \rtimes \Bigl\{
    \e^t \begin{bmatrix}
      1 & 0 \\
      t & 1
    \end{bmatrix}:t \in \R\Bigr\} 
\\
&    \Bigl\{ \begin{bmatrix}
      u & 0 \\
      0 & 0
    \end{bmatrix} : u \in \R \Bigr\} \rtimes \Bigl\{ \begin{bmatrix}
      \e^t &  0 \\
      0 & \e^t
    \end{bmatrix} : t \in \R\Bigr\} 
\\
&    \Bigl\{ \begin{bmatrix}
      0 & u \\
      u & 0
    \end{bmatrix} : u \in \R \Bigr\} \rtimes \Bigl\{ \begin{bmatrix}
      \e^t &  0     \\
      0 & \e^{t\beta}
    \end{bmatrix}:t \in \R\Bigr\} \quad \beta \in [-1,1].
  \end{align*}
           \end{enumerate}
 None of the above groups is conjugate to any other in the list modulo $\GL(2,\R)$.
 \end{Thm}


\section{The projective Lorentz representation of $\GL(2,\R)$}

From now on we assume that $d=2$, so that  $Q= \Sym(2,\R) \rtimes \GL(2,\R) $. The basic point of our construction is to realize the vector space isomorphism $\Sym(2,\R)\simeq\R^3$ as a linear isometry $\varphi:\R^3\to\Sym(2,\R)$, where $\Sym(2,\R)$ is endowed with the  inner product $ \langle \sigma , \tau \rangle = \frac{1}{2} \tr(\sigma,\tau) $. We define
 \begin{equation}
  \varphi(x,y,t) =\begin{bmatrix}
                                          t+x & y  \\
                                          y   & t-x
                                         \end{bmatrix}.
                                         \label{coords}
 \end{equation}
 Observe that $ \eta = \det\varphi = t^2 - (x^2 + y^2)$, the canonical Lorentz quadratic form. 

 We are next interested in carrying the action \ref{eq:dag_action} from $\Sym(2,\R)$ to $\R^3$.
 In other words, we are going to express \ref{eq:dag_action} in the coordinates provided by $\varphi^{-1}$.
The isometry $\varphi$ induces the map $\varphi_*:\GL(\R^3)\to\GL(\Sym(2,\R))$ given by 
$T\mapsto\varphi\circ T\circ\varphi^{-1}$  and we define the representation $\L$ as the map which makes the diagram
\[
\xymatrix{                                        &\GL(\R^3) \ar[dd]_{\varphi_*} \\
              \GL(2,\R) \ar[ur]^{\L} \ar[dr]_{\dag}   &                                \\
                                                      & \GL(\Sym(2,\R)) & } 
 \]
  commutative.
 Explicitly, we put
 \begin{equation} \label{eq:L_action}
  \L(g)  = \varphi^{-1} g^{\dag} \varphi.
 \end{equation}
Now, for $ u \in \R^3 $
 \begin{equation*}
\det(\varphi( \L(g) u)) = \det(\t g^{-1} \varphi(u) g^{-1}) \\
                        = \det(g^{-1} \ \t g^{-1}) \det(\varphi(u)) = \det(g^{-1} \ \t g^{-1}) \eta(u).
 \end{equation*}
 Hence, if $g\in\SL(2,\R)$  we have $\eta(\L(g)u)=\eta(u)$ and hence
  $\L(\SL(2,\R)) \subset \O(2,1)$. As $\SL(2,\R)$ is connected and $\L$ is continuous, the image of $\SL(2,\R)$ lies in  the connected component of the identity, that is, in the (proper, orthochronous) Lorentz group $\SO_0(2,1)$.

Next, we want to determine kernel and range of $\L$. 
 It is very easy to see that 
 \[ 
 \ker\L = \Z_2.
 \]
 Indeed, $ \ker\L = \{ h \in \GL(2,\R): \ \t h \sigma h = \sigma \text{ for all } \sigma \in \Sym(2,\R) \} $. Now,
 $ \sigma = \left[\begin{smallmatrix} 1 & 0 \\ 0 & 0 \end{smallmatrix} \right] $ yields $ h = \left[ \begin{smallmatrix} \pm1 & 0 \\ * & * \end{smallmatrix} \right] $;
 with $ \sigma = \left[\begin{smallmatrix} 0 & 0 \\ 0 & 1 \end{smallmatrix} \right] $ one obtains $ h = \left[ \begin{smallmatrix} * & * \\ 0 & \pm1 \end{smallmatrix} \right] $;
 finally, choosing $ \sigma = \left[\begin{smallmatrix} 0 & 1 \\ 1 & 0 \end{smallmatrix} \right] $ we get $ h = \pm I $.
 Conversely, it is clear that $\pm I \in \ker\L$. 
 
 As for $\Im\L $, we split $ \GL(2,\R) $ into its two connected components $\GL^+$ and $\GL^-$ and further decompose each of them, obtaining the disjoint union
\begin{equation}
\GL(2,\R) = \R^+ \SL(2,\R) \ \cup \ \R^+\!\Lambda \SL(2,\R), 
\label{CC}
\end{equation}
 where
\[
  \Lambda = \begin{bmatrix}[r]
              1 &  0 \\
              0 & -1
             \end{bmatrix}.
\]
Next we apply the Iwasawa decomposition of $\SL(2,\R)$, namely $\SL(2,\R)=NAK$.
This amounts to saying that  any element in  $\SL(2,\R)$ is the product of a shearing, a dilation and a rotation, which correspond  to the nilpotent, abelian and  compact subgroups, respectively, that is, to $N$, $A$ and $K$. The factors of a given element are unique if their order is fixed.
As it is customary, we parametrize them as
 \begin{equation} \label{SL-Iwasawa_decomposition}
  \begin{bmatrix}
   1 & 0 \\
   t & 1
  \end{bmatrix}, \qquad \begin{bmatrix}[l]
		        \e^{s/2} & 0        \\
			 0       & \e^{-s/2}
		       \end{bmatrix}, \qquad 
		       \begin{bmatrix}[r]
				             \cos(\theta/2) & -\sin(\theta/2) \\
					     \sin(\theta/2) &  \cos(\theta/2)
					    \end{bmatrix}
 \end{equation}
 with $t,s\in\R$ and $\theta\in[0,2\pi)$.
 Using \ref{eq:L_action}, one computes that their images under $\L$ are
 \begin{equation} \label{Iwasawa_decomposition}
  \begin{bmatrix}
   1 - t^2/2 & t &   t^2/2  \\
     - t     & 1 &   t      \\
     - t^2/2 & t & 1 + t^2/2
  \end{bmatrix},\quad                \ \begin{bmatrix}
                                  \cosh s & 0 & \sinh s \\
                                   0      & 1 &  0      \\
                                  \sinh s & 0 & \sinh s
                                 \end{bmatrix},\quad             \ \begin{bmatrix}[r]
                                                              \cos \theta & -\sin \theta & 0 \\
                                                              \sin \theta &  \cos \theta & 0 \\
                                                               0     &   0     & 1
                                                             \end{bmatrix}.
 \end{equation}
These, in turn, generate the nilpotent, abelian and compact Iwasawa factors of $\SO_0(2,1)$, for which we shall again use, as it is customary, the letters $N$, $A$ and $K$, respectively.
 It follows that $\L$ maps $\SL(2,\R)$ onto $\SO_0(2,1)$.
Notice also that
 $$ \L(\Lambda) = \begin{bmatrix}
                   \Lambda & 0 \\
                    0      & 1
                  \end{bmatrix}.
                  $$
 Moreover, using again \eqref{eq:L_action} one calculates that $ \L(\lambda I) = \lambda^{-2} I $ for $ \lambda \neq 0 $, hence
\[
  \L(\R^\times) = \R^+,
\]
where  $\R^\times$ denotes the non zero real  numbers.

Denote by $\O^+(2,1)$ the orthochronous subgroup of the Lorentz group, that is, the subgroup which leaves the future invariant. Its elements are the Lorentz matrices whose lower right entry is positive. It is well known that the following disjoint union holds
 \[
 \O^+(2,1) = \SO_0(2,1) \cup \L(\Lambda) \SO_0(2,1).
 \]
Using the decomposition \eqref{CC} and this equality we obtain
 \begin{equation*}
   \L(\GL(2,\R)) = \R^+ \SO_0(2,1) \cup \R^+ \L(\Lambda) \SO_0(2,1) = \R^+ \O^+(2,1),
 \end{equation*}
 which implies
\[
  \Im \L = \R^+ \! \times \O^+(2,1).
\]

Note that, if $G=\Sigma\rtimes H$ is in the class $\mathcal E_2$, then
$H$ is connected, so that $\L(H)\subset  \R^+ \SO_0(2,1) $. Hence, in the following
sections we consider only the restriction of $\L$ to the connected
component $\GL_0(2,\R)$, whose image is $\R^+ \SO_0(2,1)$.

\section{Orbits}

In order to classify the subgroups of $Q$ that belong to $\E$, the strategy is as follows. 
We fix a vector subspace $ \Sigma \subset \Sym(2,\R) $ and determine the maximal subgroup 
$H(\Sigma)$ of $\GL(2,\R) $ that leaves $\Sigma$ invariant under the action \eqref{eq:dag_action}. Then, we find all the subgroups of $H(\Sigma)$, because any such $H$ gives rise to a group $\Sigma\rtimes H$ in $\E$. For both choices, $\Sigma$ and $H$, we must of course take into account the correct notion of equivalence, namely conjugation by $MA\simeq\GL(2,\R)$. Now, 
by \ref{eq:conjugation_split}, conjugation of a group in $\E$ splits  into the action \ref{eq:dag_action} on the vector space $\Sigma$ and again ordinary conjugation on the group $H$. Thus, it is enough to select a vector space $\Sigma$ in each equivalence class modulo  the action \ref{eq:dag_action}, and then, for the same 
$\Sigma$, pick all the possible inequivalent subgroups of $H(\Sigma)$ while keeping $\Sigma$ fixed. Therefore we must  find the subgroups of $H(\Sigma)$ that are not conjugate to each other by some $g\in\GL(2,\R)$ for which $g^\dagger(\Sigma)=\Sigma$. 

The above strategy will be implemented with the aid of two observations. The first has to do with duality.
We are looking at subspaces of the three dimensional vector space $\Sym(2,\R)\simeq\R^3$.  The nontrivial ones, lines and planes, are dual to eachother under orthogonality and this duality is compatible with the semidirect structure, namely
 \begin{equation} \label{eq:dual_group}
  H(\Sigma^\perp) \ = \ \t{H}(\Sigma).
 \end{equation}
To see this notice that if $ \tau \in \Sigma^\perp $ and $ h \in H $, then for all $ \sigma \in \Sigma $
 $$ 2 \langle (\t{h})^\dag (\tau) , \sigma \rangle = \tr (h^{-1}\tau \ \t{h}^{-1}\sigma) = \tr(\t{h}^{-1}\sigma h^{-1}\tau) = 2 \langle h^\dag (\sigma) , \tau \rangle = 0. $$
 Thus, for each group $ \Sigma \rtimes H(\Sigma)\in\E $ there is a twin group $ \Sigma^\perp \rtimes \ \t{H}(\Sigma)\in\E$.
 In view of this, we focus our attention  on the case ${\rm dim}(\Sigma)=1$. This means that we are looking at   the projective space of $\Sym(2,\R)$ and we want to determine how many distinct points we need to consider, that is, how many orbits there are under the action  \ref{eq:dag_action}.
 
The second observation is that we can use geometry and carry all the above reasoning over to $\R^3$ where the action  \ref{eq:dag_action} has become the standard action of $\R^+\!\times\!\O^+(2,1)$. Hence the selection of inequivalent vector subspaces is tantamount to finding orbit representatives of the linear action of $\R^+\!\times\!\O^+(2,1)$.

 The group $\R^+\!\times\O^+(2,1)$ acts on $\R^3$ and has exactly six orbits, including the trivial one. Borrowing the terminology from
Physics, they are: the present event $(0,0,0)$, the future ($ \eta,t >
0 $), the past ($ t < 0 < \eta $), the future and the past light cone
($ \eta = 0, t > 0 $ or $ t < 0 $), the absolute elsewhere ($ \eta < 0
$).Since we are looking at   the projective spaces,
 it is enough to select a point in the upper sheet of the hyperboloid  of two sheets ($ \eta = 1, t > 0 $), one point at latitude one in the cone ($ \eta = 0, t = 1 $) and one point in the hyperboloid of one sheet ($ \eta = -1 $).
\begin{figure}[htbp]
\begin{center}
\includegraphics[width=1.8in]{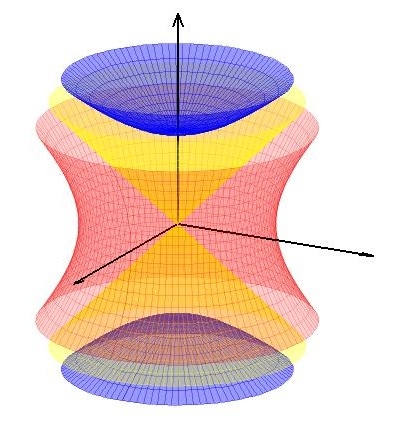}
\caption{Hyperboloids}
\label{default}
\end{center}
\end{figure}

\noindent 
 For $ u \in \R^3 \setminus 0 $ define $\hat{H}(u) $ as the projective stabilizer of $u$ with respect to the action of $ \R^+\times\!\SO_0(2,1) $, that is
 \begin{equation*}
  \begin{split}
   \hat{H}(u) = &\{ {h} \in \R^+\!\times \SO_0(2,1): \ [{h}u] = [u] \} \\
                     = &\R^+ \times \{ {h} \in \SO_0(2,1): \ [{h}u] = [u] \},
    \end{split} \end{equation*}
 where the brackets mean projective equivalence. By this definition, we have
\[
  \hat{H}(u) = \L H(\varphi(u)).
\]
 The choices that follow are convenient for computations because they take into account the Iwasawa decomposition \eqref{Iwasawa_decomposition}:
 \begin{equation*}
  \begin{array}{ccccc}
     \eta &\vline&u_\eta      \in \R^3    &\vline& \sigma_\eta \in \Sym(2,\R) \\ \hline\hline
   &\vline&&\vline&  \\
      1   &\vline& u_{1}       =  (0,0,1) &\vline& \sigma_1              =  \begin{bmatrix}
                                                                             1 & 0 \\
                                                                             0 & 1
                                                                            \end{bmatrix} \\
    &\vline&&\vline&  \\                                                                       
      0   &\vline& \ u_0         =  (1,0,1) &\vline& \sigma_0              =  \begin{bmatrix}
                                                                             1 & 0 \\
                                                                             0 & 0
                                                                            \end{bmatrix} \\
    &\vline&&\vline&  \\                                                                        
     -1   &\vline& \ u_{-1}      =  (0,1,0) &\vline& \sigma_{-1}      =       \begin{bmatrix}
                                                                             0 & 1 \\
                                                                             1 & 0
                                                                            \end{bmatrix}.
  \end{array}
 \end{equation*}
Indeed, if we denote as usual by $N$, $A$ and $K$ the Iwasawa factors of $\SO_0(2,1)$ in \ref{Iwasawa_decomposition}, then it is immediate to check that  $ K u_1 = u_1 $, $ N u_0 = u_0 $ and $ A u_{-1} = u_{-1} $.
Notice that  $\sigma_0=\varphi(u_0)/2$, but of course the constant is irrelevant.

\section{Maximal groups}

In this section we calculate the subgroups $ \hat{H}_\eta := \hat{H}(u_\eta) $ and pull them back to $\GL(2,\R)$ via  $\L$.  After that, we determine their duals.
To this end, use the Iwasawa decomposition, either in the form $NAK$ or $KNA$.
 
 If ${h}(0,0,1)$ and $(0,0,1)$ are aligned for some ${h} \in \O^+(2,1) $, then they must coincide,
 whereas ${h} (0,1,0) = (0,\lambda,0) $ implies $ \lambda = \pm1$.
 Therefore 
 $$ \hat{H}_1 = \R^+ \times \bigl(\SO_0(2,1)_{u_1}\bigr)_0, \qquad
 \hat{H}_{-1} = \R^\times \times\bigl(\SO_0(2,1)_{u_{-1}}\bigr)_0 $$
 where subscripts indicate stabilizers. Neither of the above circumstances occurs for $u_0$.
 
 Now let ${h} \in \SO_0(2,1) $, ${h} = nak $. As $K$ stabilizes $u_1$, we have ${h} u_1 = na u_1 $.
Further, $A$ moves $u_1$ along hyperbolae and $N$ along
parabolae. Hence $ na u_1 = u_1 $ implies $a=n=e$, the identity of the
group. This implies $ \SO_0(2,1)_{u_1} = K $, so that 
$\hat{H}_1 = \R^+\times K$.

 Since $A$ stabilizes $u_{-1}$, $ kna u_{-1} = kn u_{-1} $. This time $N$ rules the hyperboloid of one sheet, so $ kn u_{-1} = \pm u_{-1} $ implies $ n = e $.
 Thus $ k u_{-1} = \pm u_{-1} $ and this is only  possible if $k$ is
 either the identity or   the spatial reflection 
\[
k_\pi=\begin{bmatrix}-1 & 0 & 0 \\
                                                              0 &  -1 & 0 \\
                                                               0     &   0     & 1
                                                               \end{bmatrix}.
                                                               \]
Hence $ \SO_0(2,1)_{u_{-1}}=A\cup k_\pi A$ and $\hat{H}_{-1} = \R^+\times A$.

 Finally $ N\subset \hat{H}_0 $ because it stabilizes $u_0$, and $ A\subset\hat{H}_0 $ because  $ A u_0 = \R^+ u_0 $.
The subgroup $K$ is not in $\hat{H}_0$ since it rotates the point away from its projective orbit.
 Hence $ \SO_0(2,1)_{u_0} = NA $. Therefore
 \begin{equation*}
 \begin{split}
   &\hat{H}(u_1)   \hspace{1.4ex} = \R^+ \! \times  K \\
   &\hat{H}(u_0)   \hspace{1.4ex} = \R^+ \! \times NA  \\
   &\hat{H}(u_{-1})               = \R^+ \! \times A  .
  \end{split}
 \end{equation*}

 It is now very easy to realize these groups in
 $\GL(2,\R)_0$, because the inverse image of each
 component is well known: 
 \begin{equation*}
  \begin{split}
   H(\sigma_1)   \hspace{1.4ex} &= \R^+ \! \times \SO(2) \\
   H(\sigma_0)   \hspace{1.4ex}&= \R^+ \! \times 
   \Bigl\{ \begin{bmatrix}a&0\\b&
     a^{-1}\end{bmatrix}:a>0,b\in\R \Bigr\} 
= \Bigl\{ \begin{bmatrix}
             a & 0 \\
             b & c \end{bmatrix} : a>0,\,c> 0  \Bigr\}
\\
   H(\sigma_{-1}) &= \R^+ \! \times
   \Bigl\{\begin{bmatrix}a&0\\0&a^{-1}\end{bmatrix}
 :\,  a>0\Bigr\}
= \Bigl\{ \begin{bmatrix}
             a & 0 \\
             0 & c \end{bmatrix} : a>0,\,c> 0 \Bigr\} .
  \end{split}
 \end{equation*}
 For simplicity, we shall write $H_\eta$ for $H(\sigma_\eta)$.
 
Finally, we compute the maximal groups arising from planes. In view of \ref{eq:dual_group}, we have:
\begin{align*}
   \sigma_1^\perp&=\Bigl\{ \begin{bmatrix}u&v\\v&-u\end{bmatrix}:u, v\in\R\Bigr\} \\
   \sigma_0^\perp&=\Bigl\{ \begin{bmatrix}0&v\\v&u\end{bmatrix}:u, v\in\R\Bigr\} \\
   \sigma_{-1}^\perp&=\Bigl\{ \begin{bmatrix}u&0\\0&v\end{bmatrix}:u, v\in\R\Bigr\}.
  \end{align*}

Below is the complete list of the largest semidirect products in $ \Sym(2,\R) \rtimes \GL(2,\R) $, up to $\GL(2,\R)$-conjugations (notice that  $H_1=\!^tH_1$ and $H_{-1}=\!^tH_{-1}$):
 \begin{equation} \label{maximal_classification}
  \begin{split}
   \sigma_\eta       &\rtimes H(\sigma_\eta), \qquad \eta = 1,0,-1 \\
   \sigma_\eta^\perp &\rtimes H(\sigma_\eta), \qquad \eta = 1,-1\\
   \sigma_0^\perp &\rtimes\,\t{H(\sigma_0)}.
  \end{split}
 \end{equation}

\section{Subgroups}
 
 Our task is now to determine all the subgroups of the groups in the list \ref{maximal_classification} that belong to $\E_2$.
 Again, we can focus on the first row of the list and transpose the results.
As explained earlier, the groups we are after are the semidirect products of $\sigma_\eta$ with connected subgroups of $H_\eta$, thus we must find the latter ones. A natural approach is to look at the Lie algebra of $H_\eta$ and determine all the Lie subalgebras. In the end, exponentiation will give us the connected subgroups. Once again, the classification  is modulo conjugations that keep the normal factor $\Sigma={\rm sp}\{\sigma_\eta\}$ fixed,
namely modulo $ g \in H(\sigma_\eta) = H_\eta$.

 \subsection{$ \eta = 1 $}
The Lie algebra of $H_1$ is $\lie{h}_1 = \R + \ \lie{so}(2)$.
A natural basis is $\{I , J\}$, with $J$ as in \eqref{J}. 
The only proper subalgebras are the one-dimensional vector spaces, with trivial commutators.
  Thus, the Lie subalgebras of $\lie{h}_1$ are the points of the real projective line, whose homogeneous coordinates  with respect to the chosen basis will be denoted $[x:\alpha]$.
  Now, $\P^1$ is the disjoint union $ \R \ \cup \{\infty\} $, according to  the usual cell decomposition 
  $\{[1:\alpha]:\alpha\in\R\} \cup\{[0:1]\}$.  We must now select a single point in each equivalence class. Let us look at  the effect of conjugation by elements in $H_1$.
  Obviously $\R^\times$ gives trivial conjugations, and so does $\SO(2)$ because $ R_\theta J R_{-\theta} = J $. Further, since
\[
\Lambda R_\theta [t(I + \alpha J)] R_{-\theta} \Lambda = t(I + \alpha \Lambda J \Lambda) = t(I - \alpha J),
\]
  it is enough to consider $ \alpha \geq 0 $. 
Since
  \begin{equation*}
   \begin{split}
    &\e^{t(I + \alpha J)} = \e^t \e^{t \alpha J} = \e^{t }R_{-t \alpha} \\
    &\e^{tJ} = R_{-t}
   \end{split}
  \end{equation*}
by taking exponentials we obtain the 1-dimensional connected subgroups 
$ \{\e^{t }R_{t \alpha}:t \in \R\}$, with $\alpha \geq0$, and
$ \SO(2)$  inside $H_1$.

 \subsection{$ \eta = 0 $}
The Lie algebra $\lie{h}_0$ of $H_0$ consists of the lower triangular matrices. Pick as generators the identity matrix and 
\[
X= \begin{bmatrix}
          0 & 0 \\
          1 & 0
         \end{bmatrix} \qquad Y= \begin{bmatrix}
                                   0 & 0 \\
                                   0 & 1
                                  \end{bmatrix},
\]
which satisfy $[X,Y] = - X$. Since  $\lie{h}_0$ has dimension $3$,  Lie subalgebras of dimensions  one and two need to be considered.
  
The subalgebras of dimension one are parametrized by $\P^2$, with coordinates induced by the basis $I,X,Y$. We write
$ \P^2 = \R^2 \cup \R \cup \{\infty\} $ in terms of the standard cell decomposition 
$\{[\lambda:\mu:1]:\lambda,\mu\in\R\} \cup\{[\lambda:1:0]:\lambda\in\R\} \cup\{[1:0:0]\}$.
  Some of these points are conjugated. Indeed, take $ h \in H_0 $, hence
  $$ h = \begin{bmatrix}
           a & 0 \\
           b & c
          \end{bmatrix}, \qquad ac \neq 0. $$
  The conjugation effect on the fat cell is an affine transformation, namely
  $$ h [\lambda:\mu:1] h^{-1} = [\lambda:\frac{a}{c}\mu -\frac{b}{a}:1], $$
  so that every element in that cell can by conjugated to $[\lambda:0:1]$. Further,
  $$ h [\lambda:1:0] h^{-1} = [\frac{c}{a}\lambda:1:0] $$
  shows that every point in the second cell is conjugate either to  $[0:0:1]$ or to $[1:0:1]$.
  
Finally, the relevant exponentials are:
 \[
 \e^{tY} =  \begin{bmatrix}
                                           1 &  0  \\
                                           0 & \e^t
                                          \end{bmatrix}, \qquad
    \e^{t X} = \begin{bmatrix}
                 1 & 0 \\
                 t & 1
                \end{bmatrix}.
\]
  Therefore, up to  conjugation, the $1$-dimensional connected subgroups of $H_0$ are
  \begin{equation*}
   \begin{split}
    &\Bigl\{ \e^{t\lambda} \begin{bmatrix}
                                                    1 &  0  \\
                                                    0 & \e^t
                                                   \end{bmatrix}:t \in \R\Bigr\},\qquad\lambda \in \R \\
             & \Bigl\{ \begin{bmatrix}
                                      1 & 0 \\
                                      t & 1
                                     \end{bmatrix}:t \in \R\Bigr\}\\ 
      & \Bigl\{ \e^t \begin{bmatrix}
                                                                                                      1 & 0 \\
                                                                                                      t & 1
                                                                                                     \end{bmatrix}:t \in \R\Bigr\} \\
       & \Bigl\{e^t I:t \in \R\Bigr\}.
   \end{split}
  \end{equation*}
  
  By duality, the $2$-dimensional vector subspaces of $\lie{h}_0$ are again a copy of $\P^2$, but of course not all of them are Lie algebras.
  Let us first determine a cell decomposition.
  Orthogonal planes are easily calculated for each cell:
  \begin{equation*}
   \begin{split}
    [1:\mu:\lambda]^\perp &= \langle \mu I + X , \lambda I + Y \rangle \\
    [0:1:\lambda]^\perp   &= \langle I , \lambda X + Y \rangle \\
    [0:0:1]^\perp         &= \langle I , X \rangle,
   \end{split}
  \end{equation*}
  where the symbol on the right hand side stands for the Lie algebra generated by the indicated elements.
From
  $$ [\mu I + X , \lambda I + Y] = [X,Y] = -X $$
we infer that only the  plane corresponding to  $\mu = 0$ is a Lie algebra, namely $ \langle X , \lambda I + Y \rangle$. Every  plane in the other two cells is clearly an abelian subalgebra. Thus, we select
  $$  \langle X , \lambda I + Y \rangle, \quad \langle I , \lambda X + Y \rangle, \quad \langle I , X \rangle.$$
Next we look for possible conjugations. For $ h \in H_0 $ as before we have
  \begin{equation*}
   \begin{split}
    &h X  h^{-1} = \frac{c}{a} X \\
    &h Y h^{-1} = - \frac{b}{a}X + Y.
   \end{split}   
  \end{equation*}
 By inspection, we see that with $a=c=1$ and $b=\lambda$ we may conjugate $\lambda X+Y$ into $Y$, so that every plane in the second family is conjugate to the plane with $\lambda=0$. No other conjugation occurs. Therefore the following is a complete list of representatives
  $$  \langle X , \lambda I + Y \rangle, \quad \langle I , Y \rangle, \quad \langle I , X \rangle. $$
Observe that
 \[ \e^{t(\lambda I + Y)} = \e^{t\lambda} \begin{bmatrix}[l]
                                                                  1 &  0  \\
                                                                  0 & \e^t
                                                                 \end{bmatrix}.
                                                                 \]
We thus get the list of $2$-dimensional connected subgroups of $H_0$, up to  conjugation:
  \begin{equation*}
   \begin{split}
     & \Bigl\{ \begin{bmatrix}[l]
                                      \e^{t\lambda} &  0                \\
                                       s            & \e^{t(\lambda +1)}
                                     \end{bmatrix} :t,s \in \R\Bigr\}, \qquad \lambda \in \R\\
    & \Bigl\{ \begin{bmatrix}[l]
                                      \e^t &  0  \\
                                       0   & \e^s 
                                     \end{bmatrix}:t,s \in \R\Bigr\} \\
    &\Bigl\{ \begin{bmatrix}[l]
                                      \e^t &  0  \\
                                       s   & \e^t
                                     \end{bmatrix} :t,s \in \R\Bigr\}.
   \end{split}
  \end{equation*}
  
 \subsection{$ \eta = -1 $} 
The Lie algebra $\lie{h}_{-1}$ of $H_{-1}$ is composed by all diagonal matrices, a natural basis of which is given by
  $$ Z = \begin{bmatrix}
           1 & 0 \\
           0 & 0
          \end{bmatrix}, \qquad Y = \begin{bmatrix}
                                    0 & 0 \\
                                    0 & 1
                                   \end{bmatrix}. $$
The proper subalgebras of $\lie{h}_{-1}$ are the points of $\P^1$. We decompose again with respect to the basis as $\{[1:\beta]:\beta\in\R\} \cup\{[0:1]\} $.  Conjugations by invertible diagonal matrices are trivial, while $R_{\pi/2}$ exchanges coordinates:
\[
    R_{\pi/2}[1:\beta] R_{-\pi/2} = [\beta:1], \qquad
    R_{\pi/2}[0:1]R_{-\pi/2} = [1:0].
\]
  Hence $[1:\beta]$, with $\beta\in[-1,1]$, is enough.  It is easy to see that
  \[ 
  \exp(t[1:\beta])=\e^t\!Z + \e^{t\beta}\!Y,
  \] 
  whence we obtain all the connected 1-dimensional subgroups of $H_{-1}$, up to conjugation
  \begin{equation*}
   \Bigl\{ \begin{bmatrix}
                               \e^t &  0     \\
                                0   & \e^{t\beta}
                              \end{bmatrix}:t \in \R\Bigr\}, \qquad \beta \in [-1,1].
  \end{equation*}

\end{document}